\documentclass[12pt]{amsart}
 \usepackage[francais]{babel}
\usepackage[applemac]{inputenc}
\usepackage{graphicx}

\DeclareMathOperator{\Card}{Card}

\newcommand{\zd}{\ensuremath{\mathbf Z_{12}}}

\newcommand{\zc}{\ensuremath{\mathbf Z_{c}}}
\newcommand{\N}{\ensuremath{\mathbf N}}
\newcommand{\Z}{\ensuremath{\mathbf Z}}
\newcommand{\C}{\ensuremath{\mathbf C}}

 \let\mc=\mathcal
 \newtheorem{Thm}{Théorème}
 \newtheorem{Def}{Définition}
  \newtheorem{Prop}{Proposition}
  \newtheorem{Lemme}{Lemme}
  \newtheorem{Rem}{Remarque}

\title{Gammes Bien R\'eparties et Transform\'ee de Fourier discr\`ete}
\author{Emmanuel Amiot \\
St Nazaire, France}

\begin{document}
\maketitle
{\sl Un des concepts les plus séduisants avancés ces dernières années
par les théoriciens américains de
la musique est celui des gammes bien  réparties (\og Maximally Even Sets\fg),
qui étend et précise celui des gammes monogènes (\og Generated Scales\fg)
ou au concept moins connu mais fécond de gamme bien formée (\og Well Formed\fg).
L'objectif essentiel de cet article est de présenter ces notions au public francophone.
Sa principale originalité consiste à partir, contrairement à l'école américaine,
d'une définition purement (géo)métrique en termes de transformée de Fourier,
sans référence à la notion ordonnée de gamme musicale. Il peut paraître ironique
que la transformée de Fourier discrète resurgisse dans ce contexte plutôt algébrique,
voire ensembliste, alors que l'étude des spectres (continus) est depuis longtemps un
cheval de bataille en théorie du son !
}\medskip

{ Notations : 
\begin{enumerate}
\item
$\lfloor x\rfloor$ dénote la partie entière du réel $x$, 
$\lceil x\rceil$ est l'entier immédiatement supérieur.
\item
$\mathbf Z$ est l'ensemble des entiers relatifs, $\zc$ le groupe cyclique
à $c$ éléments.
\item
$c\wedge d$ est le plus grand commun diviseur de $c$ et $d$.
\item
$d\mid c$ signifie que $d$ divise $c$.
\end{enumerate}
}

\section{Que veut-on dire par \og bien répartie\fg ?}

Dans le langage de Claude Debussy comme de bien d'autres compositeurs,
on trouve aussi bien la gamme traditionnelle, majeure, que la gamme par tons
ou la gamme pentatonique; dans plusieurs pièces (\textsl{L'Isle Joyeuse, Voiles…})
on trouve plusieurs de ces modes en concurrence, ou en complément. Le présent article
présente, sous un angle inédit, une propriété qui est commune à ces trois gammes
et les caractérise sous certaines conditions.

Nous travaillons dans le contexte d'un \textbf{total chromatique}
modélisé par le groupe cyclique $\zc$\footnote{ce qui modélise le fait
psychoacoustique que le cerveau tend à reconnaître comme \og égales\fg\ des
notes distantes d'une octave, i.e. dont le rapport des fréquences vaut 2.}  
à $c$ éléments, qui s'interprète musicalement comme une division
de l'octave en $c$ parties égales, ou encore aux rapports de fréquence
de la forme $2^{k/c}, k\in\N$. Il existe en musique
des tempéraments inégaux, bien sûr, et l'on peut
étendre au moins une partie de la théorie qui suit à de tels objets, avec des résultats 
intéressants, mais ce n'est pas ici notre propos.
Une \textbf{gamme} sera tout simplement un sous-ensemble de cardinal $d$ 
de $\zc$ (ces notations sont celles notamment de
\cite{Quinn}). Notre souci est de choisir ces $d$ notes parmi $c$ de telle sorte 
qu'elles soient les plus également réparties possible. Les pionniers de cette notion ont 
découvert avec délice que l'on retrouve ainsi des gammes très particulières, comme la gamme 
majeure ou la pentatonique (gamme chinoise) qui trouvent ainsi une légitimité supplémentaire.

\subsection{Bien poser la question}\ 
Il convient de préciser ce que l'on entend par \og le mieux réparti possible\fg.

Ce n'est pas une question si évidente, comme évoqué par \cite{BD}. Ceci
peut sembler étrange, car l' idée en
est intuitive et suffisamment simple pour être appréhendée par des lycéens (\cite{JD}).

Mais si l'on propose, par exemple, la définition naturelle consistant à maximiser la somme des
 intervalles entre tous les éléments (l'intervalle entre $a,b\in\zc$ étant 
 $\max\limits_{k\in\Z} |a-b+k\, c|$, la plus
 petite valeur de $|a-b|$ quand on fait varier le choix des représentants des classes
 $a, b$ modulo $c$), on obtient des gammes tout à fait contraires à l'intuition.
Ainsi par exemple : parmi les pentacordes (gammes de cinq notes parmi 12), 
les diverses valeurs possibles de cette somme sont
40, 48, 52, 56, 60, 64, 68 et 72. Si la valeur minimum de 40 correspond, conformément
à l'intuition, aux pentacordes chromatiques (comme Do Do$\sharp$ Ré Ré$\sharp$ Mi),
le maximum de 72 est atteint pour de très nombreuses (21 !) gammes différentes,
comme (0 1 2 6 7) ou  Do Do$\sharp$ Ré Fa$\sharp$ Sol - il n'y a donc pas unicité
de la solution, qui ne caractérise donc pas la gamme pentatonique\footnote{ Par
rapport au résultat général de \cite{Ising} mentionné plus loin, cela tient à ce que la fonction utilisée
d'évaluation n'est pas \textbf{strictement } convexe.}.

La définition usuelle par l'école américaine présuppose que l'on range les éléments
de la gamme, pour pouvoir calculer les intervalles entre une note et ses voisines immédiates. 
Elle porte donc une forte connotation culturelle ($\zc$ n'a pas une vraie relation d'ordre !
on se souvient ici \textbf{inconsciemment} que le modèle de $\zc$
est une réduction du domaine totalement ordonné des hauteurs),
et nous ferons ici le choix de proposer une définition purement abstraite, qui
s'applique à un \textbf{sous-ensemble} à $d$ éléments et non à un $d-$uplet.

Une méthode intuitive avec \og des points blancs et des points noirs\fg\ (\cite{CD}),
bien qu'elle consiste à étudier justement l'ordre dans lequel se retrouvent les points,
mérite une mention spéciale : on considère un polygone régulier à $d$ sommets blancs et un autre
avec $c-d$ sommets noirs. Alors on superpose les deux, convenablement tournés pour que deux sommets
ne puissent coïncider; l'ordre alors obtenu donne une gamme bien répartie\footnote{
Aussi bien pour le sous-ensemble noir que pour le blanc.}, que l'on peut expliciter
en ajustant les positions des points sans changer leur répartition. Clough et Douthett montrent que
ce procédé équivaut à leur définition. 


En tant que mathématicien, on peut rechercher une définition moins \og bricolée\fg.
À titre de motivation citons un résultat  relativement récent (\cite{Toth})

\begin{Thm}[(Toth, 1956)]
   On considère $d$ points sur le cercle unité; la somme des distances \textbf{euclidiennes}
   (i.e. en ligne droite) entre tous ces points est maximale si, et seulement si, la figure qu'ils
   forment est un polygone régulier.
\end{Thm}
\begin{proof}
  Elle est non triviale. Observons que ce résultat est faux en dimension supérieure
  (c'est encore un problème ouvert en général que de trouver $d$ points sur la sphère $S^2$
  maximisant la somme de leurs distances mutuelles). En revanche il est facile
  (calcul différentiel élémentaire) de montrer que les polygones réguliers sont
  ceux qui maximalisent le \textbf{périmètre}.
\end{proof}

Nous pourrions, en nous inspirant de ce lemme, chercher les sous-ensembles à $d$
éléments de $\zc$, identifié au $c-$polygone régulier, qui
maximisent la somme des distances
euclidiennes, ce qui est une façon de chercher à approximer un polygone régulier
à $d$ sommets. Ce serait une bonne idée dans la mesure où l'on retrouverait
ainsi les sous-ensembles particuliers appréciés des musiciens (\cite{Ising}),
et une unique solution pour chaque valeur de $d$ (à symétries près). 
Si l'on reprend le cas des pentacordes dans l' univers des douze sons, la
gamme pentatonique est maximale pour ce critère, avec une somme égale à
{$\oldstylenums{30.5758}$}, et c'est l'unique solution à translation près.

Mais ces distances euclidiennes, ces longueurs de cordes, n'ont pas de signification 
musicale\footnote{ On prend le double du sinus de l'intervalle, avec un facteur d'échelle\dots} 
et font appel à la structure euclidienne du plan dans lequel on dessine la figure;
nous préfèrerons donc une définition sans doute plus technique, mais à notre avis plus pure,
qui va dans le même sens, à savoir chercher à approximer un $d-$polygone régulier 
par une sous-figure du $c-$polygone régulier. Ceci est plus proche du procédé avec les
points noirs et blancs.
Mais dans ce but nous allons plutôt faire un retour aux sources de l'étude des propriétés
internes des gammes (ou accords) telle qu'initiée par \cite{Lewin59}.

\subsection{Transformée de Fourier}\ 

L'idée géniale d'exploiter la transformée de Fourier discrète pour inventorier
les intervalles dans un accord (ou une gamme) remonte à David \textsc{Lewin} en 1959.
et est à l'origine de nombreux concepts qui ont marqué la recherche américaine depuis
1959. Nous choisissons la même définition que lui parmi toutes les définitions possibles:

\begin{Def}
   La transformée de Fourier de $f:\zc\to \mathbf C$ est
  $$\mc{F}(f):t\mapsto  \sum\limits_{k\in\zc} f(k) e^{-2 i k \pi t/c}$$
  Plus particulièrement, la transformée de Fourier de $A\subset \zc$ 
  sera la transformée de Fourier de la fonction caractéristique $1_{A}$
  du sous-ensemble $A$\footnote{Si l'on s'intéresse aux fonctions du cercle
  $S^1$ à valeurs dans $\C$, on peut voir cela comme la transformée de
  \textsc{Fourier} d'une distribution, à savoir
  un peigne de \textsc{Dirac} $\sum_{k\in A} \delta_{k}$.}:
   $$
     \mc{F}_{A}:t\mapsto \sum_{k\in A} e^{-  2i\pi k t /c}
   $$
\end{Def}
Quelques exemples:
\begin{enumerate}
  \item
  $ \mc{F}_{\zc}$, la transformée de Fourier de toute la gamme chromatique, est
  $\sum\limits_{k=0}^{d-1} e^{- 2 i\pi k t  /c}
         =\dfrac{1 -  e^{- 2 i \pi t  }}{1 -  e^{- 2 i\pi t /c}}$.
  Cette fonction est nulle sur $\zc$ sauf quand $t=0$.
  On voit bien sur ce calcul que seule compte la classe de l'indice $k$ modulo $d$,
  ce qui est adéquat.
  \item
  Considérons la septième diminuée D7 = (0 3 6 9) dans l'univers à 12 notes. Alors
  $$\mc{F}_{D7}:t\mapsto \sum_{k=0}^{3} e^{-(3 k) 2i \pi t / 12}=
  \sum_{k=0}^{3} e^{-k i\pi t / 2}=
  \dfrac{1 - e^{- 2 i \pi t}}{1 - e^{- i \pi t/2}}$$
  On remarque que $\mc{F}_{D7}(4)=4$. Observons aussi que $A$ est $3-$périodique:
  le principe même de la transformée de Fourier est que \textbf{le
  coefficient de Fourier $\mc{F}(f)(k)$ mesure à quel point
  la fonction $f$ est $c/k - $périodique}. Ceci résulte du principe de la transformée de Fourier
  inverse (théorème de \textsc{Plancherel}): on a toujours
  $$
     \forall t\in\zc \qquad
     f(t) = \frac1c \sum_{k=1}^c \mc{F}(f)(k) e^{+ 2i\pi k t /c}
  $$
  et donc plus $\mc{F}(f)(k)$ est gros, plus $f$ est \og concentrée\fg\
   autour de $t\mapsto e^{+ 2i\pi k t /c}$.
  \item
  Le \textbf{module} de la transformée de Fourier est invariant par translation
  et inversion :
  $$
  \mc{F}_{A-p}(t) = e^{2 i p t/c} \mc{F}_{A}(t)
     \qquad 
  \mc{F}_{A}(-t) = \overline{\mc{F}_{A}(t)}
  $$
  Cette propriété est importante pour les musiciens, dont les notions sont invariantes
  par transposition (= translation dans le domaine des hauteurs) et parfois par
  inversion (où la quinte devient une quarte).
  \item
  Pour toute partie $A\in\zc$ on a  $\mc{F}_{A}(0) = \Card\, A$.

\end{enumerate}
La puissance de cet outil est considérable, on la voit sur les propriétés suivantes:

\subsection{Propriétés des transformées de Fourier de parties de $\zc$}
\begin{Thm} \label{complementaire}
  La transformée de Fourier d'une partie
  et celle de son complémentaire sont opposées
  en $k=1… d-1$:
  $$
  \mc{F}_{A}(k) + \mc{F}_{\zc\setminus A}(k) = 0\quad
  \forall k=1…c-1
  $$
\end{Thm}
\begin{proof}
Cela provient, par linéarité, de $ \mc{F}_{A} + \mc{F}_{\zc\setminus A} = \mc{F}_{\zc}$ 
et de ce que (cf. supra)
$\mc{F}_{\zc}(k) = 0 $ pour $ k=1… c-1$.
\end{proof}

Observons que la remarque faite ci-dessus sur la septième diminuée se généralise:

\begin{Prop}
    Considérons une division de $\zc$ en $d$ parties égales:
    $A = \{0, m=\dfrac{c}d, 2m, … (d-1)m\}$. Alors la transformée de Fourier
    de $A$ est maximale au point $t=d$. Plus précisément,
    $$
    \mc{F}_{A}(t) = 
       \begin{cases}
            d &\text{ si } d\mid t\\
            0 &\text{ sinon }
       \end{cases}
    $$     
\end{Prop}
\begin{proof}
  En effet, on a 
  $$\mc{F}_{A}(t) = \sum_{k=0}^{d-1} e^{-2 i k \pi t  m/c}
  = \sum_{k=0}^{d-1} e^{-2 i k \pi t / d} 
  = \dfrac{1 - e^{-2 i \pi t}}{1 - e^{-2 i \pi t/d}}$$
  pour $t$ non multiple de $d$. En revanche pour $t=d$ (et multiples) on a
  $\mc{F}_{A}(t) = 1+1+\dots 1 = d$.
\end{proof}
Il est notable que la division régulière de $\zc$ donne une valeur absolue maximale de
la transformée de Fourier \textbf{parmi tous les sous-ensembles de cardinal $d$ donné}.
Nous allons toutefois donner une raison supplémentaire
de nous intéresser à cette transformée, qui est qu'elle permet de mesurer l'importance de la 
présence entre les éléments de $A$ d'intervalles proches de $c/d$.

\begin{Thm}[(Contenu Intervallique)]\ 

  Définissons le \textbf{contenu intervallique} de $A\subset \zc$ 
  par le nombre d'occurences de chaque intervalle entre deux éléments de $A$, 
  ce qui se résume par la fonction:
  $$
    IC_{A}(k) = \Card \{ (x,y)\in A\times A \mid x+k=y\}
  $$
  Alors la transformée de Fourier de la fonction $IC_{A}$ est
  ${|\mc{F}_{A}|}^2$.
\end{Thm}
C'est la remarque fondamentale de \textsc{David Lewin} dans
son tout premier article (\cite{Lewin59}), remarque
qu'il a énoncée (en deux lignes seulement) dans le contexte plus général des
\textbf{relations intervalliques} entre deux parties. 

\begin{proof}
   IC est un produit de convolution (que nous noterons  $\star$) 
   entre la fonction caractéristique de $A$ et celle de $-A$:
    $$
    IC_{A}(k) = \Card \{ (x,y)\in A\times A \mid x+k=y\}
    = \sum_{\substack{y\in A\\y-k\in A}} 1
    = \sum_{y\in\zc} 1_{A}(y)\times 1_{-A}(k-y)
    = \bigl(1_{A}\star 1_{-A}\bigr)(k)
  $$
  Il est bien connu que la transformée de Fourier d'un produit de convolution
  est le produit ordinaire des transformées de Fourier :
  $\mc{F}(f \star g) = \mc{F}(f)\times \mc{F}(g)$.
  Dans $\zc$, cela dérive du calcul suivant:
  \begin{multline*}
    \mc{F}(f \star g)(k) = \sum_{k\in\zc} (f\star g)(k) e^{-2 i k \pi t/c}
    = \sum_{k\in\zc}\sum_{\ell\in\zc} f(\ell) g(k-\ell) e^{-2 i k \pi t/c}
    = \sum_{\ell\in\zc}\sum_{k-\ell\in\zc} f(\ell) g(k-\ell) 
                e^{-2 i (k-\ell) \pi t/c} e^{-2 i \ell \pi t/c}\\
                    =  \sum_{\ell\in\zc}f(\ell)e^{-2 i \ell \pi t/c} \times
                \sum_{k-\ell\in\zc}  g(k-\ell) e^{-2 i (k-\ell) \pi t/c} 
                                 = \mc{F}(f)(k)\times \mc{F}(g)(k)
  \end{multline*}
  Appliquant cela à $f = 1_{A}$ et $g = 1_{-A}$ on obtient le résultat 
  (avec $\mc{F}_{-A} = \overline{\mc{F}_{A}}$).
\end{proof}
On déduit au passage de ce résultat une très jolie preuve du théorème de l' hexacorde 
de Milton \textsc{Babbitt}:

\begin{Thm}[de l'hexacorde (Babbitt, 1964)]\ 

   Deux parties complémentaires et de même cardinal de $\zc$ (pour $c$ pair)
   ont même contenu intervallique.
\end{Thm}
En effet, si l'on considère deux parties complémentaires et de même cardinal, 
alors leurs transformées de Fourier sont opposées 
(pour $t=1…c-1$) ou égales (en $t=0$), en particulier les transformées de Fourier
de leurs IC sont égales, et donc les IC aussi
(par transformée de Fourier inverse).\medskip

Nous avons donné ce théorème sur le contenu intervallique parce qu'il
 montre bien l'importance de la transformée de Fourier
quant aux intervalles présents dans une partie de $\zc$:
la valeur de la transformée de Fourier de $IC_{A}(t)$ pour $t=m$ dit
dans quelle mesure cette fonction IC est $c/m-$ périodique: 
considérons le cas extrême où $\mc{F}(IC_{A})$ serait nulle, sauf
en $d$ (et multiples), où elle vaudrait $d$: alors d'après ce qui précède, $A$ 
doit diviser régulièrement $\zc$.

Tout ceci nous convainc que la transformée de Fourier, et plus particulièrement
son module, donne une bonne idée d'à quel point une partie $A\subset \zc$ est proche
d'un polygone régulier. Cette philosophie de repérer les divisions \og archétypales\fg\
du total chromatique est celle de la première partie de la remarquable 
dissertation de Ian Quinn (\ref{Quinn}) qui a motivé le présent exposé.

\subsection{Enfin une définition}
\begin{Def}
  Dorénavant nous dirons que $A\in\zc$, de cardinal $d$, est
  \textbf{bien répartie} si la quantité
  $$\sqrt{\mc{F}(IC_{A})(d)} = |\mc{F}_{A}(d)|$$ 
  est maximale parmi toutes les parties de cardinal $d$:
  $$
     |\mc{F}_{A}(d)| ≥  |\mc{F}_{A'}(d)| \quad
     \forall A'\subset\zc, \Card\, A'=d
  $$
\end{Def}
Observons que cette définition (contrairement à celle que donnent \cite{CD} \textsl{et alia})
porte sur \textbf{l'ensemble $A$}, sans considérer du tout les éléments de $A$ comme
ordonnés. Une partie à $d$ éléments qui réalise ce minimum
(il en existe car il n'y a qu'un nombre fini de parties à $d$ éléments parmi $c$)
sera dite \textbf{gamme bien répartie} ou GBR.

Résumons ici quelques remarques (démontrées \textsl{supra})
qui motivent cette définition: 
\begin{itemize}
     \item
     $|\mc{F}_{A}| = |\mc{F}_{A+\alpha}|=|\mc{F}_{-A}|$ pour tout $\alpha$ 
     (invariance par translation et inversion);
     \item
     $|\mc{F}_{A}(t)| ≤ d=\Card\, A$, pour tout $t\in\zc$ et même dans $\mathbf R$.
     \item
      $|\mc{F}_{A}(p)| = d$ si et seulement si toutes les exponentielles $e^{-2i k\pi t/c}$ 
      qui figurent dans $\mc{F}_{A}$ ont même argument (modulo $2\pi$), 
      i.e. pour $a,b\in A$ la quantité $(b-a)p/c$ est toujours un entier.
      Ce qui prouve la propriété suivante, réciproque de la proposition précédente:
\end{itemize}
\begin{Thm}
    Soit $A\subset \zc$ de cardinal $d$; si $|\mc{F}_{A}(d)| = d$, alors
    $d$ divise $c$, et $A$ est un polygone régulier (i.e.
    $A = \{a_{0}, a_{0}+c/d,\dots a_{0}+(d-1)c/d\}$).
\end{Thm}
\begin{proof}
   Fixons une origine $a_{0}\in A$: on a dit que tous les $(a - a_{0})d/c, a\in A$ prennent 
   $d$ différentes valeurs \textbf{entières} (modulo $c$), donc le pgcd $g$ des
    $a-a_{0}$ est multiple de $c/d$: $c/d \mid g$ et $g≥c/d$.
  
  Par ailleurs, il y a $d$ multiples de $g$, compris entre 0 et $c-1$,
  qui sont distincts, 
  à savoir les $a-a_{0}, a\in A$.\footnote{On choisit le
  représentant de la classe modulo $c$ qui est compris entre 0 et $c-1$} 
  Donc $g$ n'excède pas $c/d$.
  D'où $g=c/d$ et les $a-a_{0}$ sont les multiples de $c/d$.   
\end{proof}

Une autre propriété tombe en une ligne :

\begin{Thm}
   Le complémentaire d'une gamme bien répartie est aussi bien répartie.
\end{Thm}
\begin{proof}
Soient $A$ et $B$ deux parties complémentaires (non vides) de $\zc$, de cardinaux 
$d$ et $c-d$ respectivement. Alors d'après le théorème \ref{complementaire}
    $$
      |\mc{F_{A}}(d) | = |-\mc{F_{B}}(d) | = |-\overline{\mc{F_{B}}(-d)} |
      = |-\overline{\mc{F_{B}}(c-d)} |= |\mc{F_{B}}(c-d) |
    $$
 et l'un est maximal (parmi toutes les parties à $d$ éléments) si et seulement si
 l'autre l'est (parmi toutes les parties à $(c-d)-$ éléments).
\end{proof}
En conséquence nous supposerons dorénavant que $c>d/2$ (le cas $c=2d$ étant trivial).

Il est aussi évident que le rétrograde $-A$ d'une partie bien répartie $A$ l'est aussi,
de même pour les translatés $A+p$.
Une propriété n'est toujours pas évidente à partir de cette définition, c'est l'unicité (à
translation près) d'une partie bien répartie de cardinal donné. 
Cela résultera des propriétés données dans la partie suivante.

\section{La formule génératrice}
\subsection{Quelques lemmes}\
  
Partant de notre définition, nous voulons retrouver qu'une partie bien répartie 
est donnée (cf. \cite{CM}, généralisé par \cite{CD}) par une formule du type
$$
    a_{k} = J_{c,d}^\alpha\lfloor\dfrac{k c + \alpha}d \rfloor
$$ 
où $a_{k}$ dénote le $k\raise2pt\hbox{ème}$ élément de $A$
et $\lfloor x \rfloor$ est la partie entière de $x$.
  
Nous avons besoin d'un lemme géométrique.
Nous voulons exprimer plus précisément que, par définition,
tous les $e^{i k d 2\pi/c}$ doivent être aussi proches que possible, puisque l'on
veut en maximiser la somme; de fait,

\begin{Thm}
   Considérons un ensemble $A$ de $d$ points \textbf{distincts} dans $\zc$; la somme
   $|\sum\limits_{k} e^{i a_{k} 2\pi/c}|$ est maximale 
   quand les points sont consécutifs dans \zc, i.e.
   $A$ est un translaté de $\{0,1,2,…, d-1\}$
   (on exige que les points soient \textbf{serrés} les uns contre les autres).
\end{Thm}

Cette propriété pourrait être adoptée comme définition des gammes bien réparties,
elle a bien sûr été remarquée comme corollaire des définitions usuelles (\cite{CD}).
L'implication présentée ici est plus originale, elle repose\footnote{
 Ce procédé d'amélioration par itération est utilisé, dans un contexte plus compliqué,
 pour démontrer le premier résultat de \cite{Ising}.} sur le

\begin{Lemme}[des points bien serr\'es]\

   On prend $d$ points $a_{1}, … a_{d}$ sur le cercle unité, et on déplace
   l'un de ces points
   (disons $a_{1}$) vers la somme de tous les $a_{i}$, c'est à dire que l'on remplace
   $a_{i}$ par un $a'_{1}$ situé sur l'arc de cercle (minimal)
   entre $a_{1}$ et $\sum\limits_{k=1}^d a_{k}$.
   Alors $|a'_{1}+a_{2}+…a_{d}| ≥|a_{1}+a_{2}+…a_{d}| $.

\end{Lemme}
\begin{proof}
La somme augmente car l'angle entre $\sum\limits_{k=1}^d a_{k}$
et $a'_{1} - a_{1}$ est aigu: en effet, si on pose que l'argument de $\sum a_{k}$ est 0
et que l'on prend $a_{1}, a'_{1}$ \og au dessus\fg, avec l'orientation usuelle du cercle; 
alors par convention il vient
$$ 0 < \arg a'_{1}=\theta'<\theta= \arg a_{1}≤\pi
\qquad \arg (a'_{1} - a_{1}) = \arg (e^{i\theta'}- e^{i\theta'}) = 
\dfrac{\theta+\theta'-\pi}2\in ]-\frac\pi2, 0[$$
On peut aussi démontrer ce lemme en projetant sur la droite dirigée par $\sum a_{k}$,
ce qui donne la somme des cosinus des angles avec cette direction,
et arguer que la fonction cosinus est décroissante sur $[0,\pi]$. 
Donc la somme des cosinus va augmenter quand on remplace $a_{1}$ par $a'_{1}$,
Or le module de la somme $a'_{1}+a_{2}+…$ est encore supérieur (la direction
aura changé) à cette somme des cosinus, donc augmente \textsl{a fortiori}. 
\end{proof}


À partir de ce lemme, la démonstration du théorème est algorithmique:
partant d'une configuration des $e^{i a_{k} 2\pi/c}$ qui comporte des trous
($A$ n'est pas translatée de $\{0,1, …, d-1$), on effectue un déplacement comme 
dans le lemme, 
bouchant un trou avec le point voisin (le voisin vers l'extérieur, i.e.
dans la direction opposée par rapport à la somme totale). 
On peut itérer cela jusqu'à ce qu'il n'y ait plus de trous, i.e. jusqu'à ce que
tous les points soient consécutifs. 
Cela se fait en un nombre fini d'étapes (le module de la somme augmente
strictement à chaque fois, et ne peut prendre qu'un nombre fini de valeurs).
D'où le théorème. \qed

Ceci nous permet de règler le cas où $c$ et $d$ sont premiers entre eux, 
cas musicalement capital pour les musiciens (gammes à 5 ou 7 notes parmi 12, par exemple !).

\begin{Lemme}
    Supposons $d\wedge c=1$ et soit $A = \{ a_{1}, … a_{d}\}$ un ensemble
    bien réparti.
    Alors il existe un décalage $k$ tel que $d\, A$ s'écrit $k+ \{1, 2, … d\}$ et il y a
    exactement $c$ tels ensembles $A$ de cardinal $d$, tous déduits l'un de l'autre 
    par translation.
\end{Lemme}

Cette propriété a déjà été observée, notamment par divers chercheurs italiens
(\cite{Cafagna}) mais elle n'a pas été utilisée à ma connaissance comme 
\textbf{caractérisation} des gammes bien réparties.
Comme on l'établira ci-dessous, à un détail technique près cela fait aussi
l'affaire même si $d\wedge c>1$.

\begin{proof}
 Le seul point à préciser avant d'appliquer le Lemme ci-dessus 
 aux points $e^{2i \pi a_{k} d/c}, a_{k}\in A$, est que
 les éléments de $d\, A$ sont bien distincts: ceci résulte de ce que
 $x\mapsto d\times x$ est bijective dans $\zc$; vu que $d\wedge c=1$.
 Le lemme nous dit que les éléments de $d\, A$ doivent être consécutifs, i.e.
$d\, A$ est translaté de  $\{1, 2, … d\}\subset \zc$ et donc
$A$ est translaté de $d^{-1}\times\{1, 2, … d\}$ où $d^{-1}\in\subset \zc$.
\end{proof}

\subsection{Engendrer toutes les gammes bien réparties}\

\subsubsection{Cas $d\wedge c = 1$}
\begin{Thm}\label{premier}
   Avec les hypothèses précédentes ($d\wedge c = 1$), une gamme
   bien répartie $A$ est une \textbf{gamme monogène} (\og Generated Scale\fg):
   à une translation près, 
   $$
     \exists f\in\zc \mid\, A = \{ f, 2 f, … d\,f\}
   $$
\end{Thm}
\begin{proof}
   Avec les notations précédentes, on prend $f=d^{-1}\in\zc$, puis
   $A = f\times\{1, 2, … d\} $.
 \end{proof}
Ainsi par exemple les gammes majeure et pentatonique, qui sont les gammes bien
réparties de cardinal 7 et 5 respectivement dans \zd, forment une portion du \textbf{cycle
des quintes}: elles sont engendrées par $f=7$, car
$$
  \{0,2,4,5,7,9,11\} = \{-7, 0,7,14,21,28,35\} \pmod{12}
$$   
La relation $f=d^{-1}$ explique bien pourquoi ceci ne fonctionne
que quand $d\wedge c=1$. Le fait que les 5 premières notes de la gamme majeure constitue
une gamme pentatonique est remarquable, et tient à ce que $f$ (ici les quintes) comme $-f$
(ici les quartes) génèrent aussi bien les mêmes GBR, puisque comme on l'a vu 
le rétrogradé d'une GBR est encore une GBR.
Plus généralement on a l'énoncé

\begin{Prop} \label{inclusComplement}
   Soit $d<c/2, d\wedge c=1$; alors une GBR de cardinal $d$ parmi $c$ est incluse dans une 
   (plus précisément dans $c-2d+1$ distinctes) GBR de cardinal $c-d$,
   c'est à dire qu'à translation près une GBR est incluse dans sa complémentaire.
\end{Prop}
Une célèbre illustration de cette propriété (qui semble curieusement inédite alors qu'elle est
en filigrane dans \cite{CD}) est l'étude 5 en sol bémol majeur, opus 10, de Chopin, dans 
laquelle la main droite ne joue que des touches noires, la main gauche établissant dans le même
temps diverses tonalités compatibles:


Enfin, cette propriété de monogénéité, avec $f=-d^{-1}$ (le signe servant à préserver le même ordre 
sur le cercle trigonométrique entre les $1,2…d$ et les éléments de $A$)
a été mise en avant par \cite{Cafagna} \textsl{et alii} en tant qu' \textbf{index},
au sens topologique (Poincaré). 
Avec notre définition métrique des gammes bien réparties, cette propriété est 
encore plus riche de sens.\medskip

Il est maintenant possible de prouver la formule de \cite{CD}, dans
ce cas $d\wedge c = 1$.

\begin{Thm}
    Nous supposons toujours que $d\wedge c=1$ et que
     $A = \{ a_{1}, … a_{d}\}$ est une gamme bien répartie.
    Alors à translation près, $A$ est l'ensemble 
    $-\{0, \lfloor \dfrac{c}d\rfloor,… \lfloor \dfrac{(d-1)c}d\rfloor\}$.
\end{Thm}
\begin{proof}
   Nous savons que $d\, A = \{0,… d-1\} \in\zc$ (à translation près).
   
   Considérons les nombres rationnels suivants 
   $0, -\dfrac{c}d,… -\dfrac{(d-1)c}d$:
   leurs parties entières sont distinctes, car $c>2d$; de même quant à 
   leurs parties fractionnaires, cela parce que $d$ est premier avec $c$. 
   
   Les résidus modulo $c$ des quantités $d\lfloor \dfrac{k c}d\rfloor$, identifiés
   à leur représentant dans $[1-c, 0]$, tombent entre $-d+1$ et 0, car 
   $$\dfrac{k c}d - 1 < \lfloor \dfrac{k c}d\rfloor≤ \dfrac{k c}d 
   \qquad
   k\,c - d< d  \lfloor \dfrac{k c}d\rfloor≤  k\,c 
   $$
   Quand $k$ varie de 1 à $d$, les $d$ entiers $k\,c - d\lfloor \dfrac{k c}d\rfloor$
   se retrouvent donc entre 0 et $d-1$; de plus ils sont distincts, puisque comme évoqué ci-avant
   les parties fractionnaires des $\dfrac{k c}d$ sont distinctes: 
   $$
     k\,c - d\lfloor \dfrac{k c}d\rfloor =    d\bigl(\dfrac{k\,c}d - 
   \lfloor \dfrac{k c}d\rfloor \bigr) =  d\cdot \text{frac}\left (\dfrac{k c}d\right)
   $$
   donc ces $d$ entiers forment exactement l'ensemble $0,… d-1$, c'est à dire que,
   modulo $c$,
$$
      d\times \{  0, -\lfloor \dfrac{c}d\rfloor,… -\lfloor\dfrac{(d-1)c}d\rfloor\} 
    \equiv \{0,… d-1\}  \quad\text{et donc}\quad   
     A = d^{-1}\times \{0,… d-1\} 
     = -\{0, \lfloor \dfrac{c}d\rfloor,… \lfloor\dfrac{(d-1)c}d\rfloor\}
$$
 C'est bien la formule annoncée.
\end{proof}
\begin{Rem}
 Considérons une famille similaire, mais sans le signe - :
 $$\{0, \lfloor \dfrac{c}d\rfloor,… \lfloor\dfrac{(d-1)c}d\rfloor\}$$
 Comme elle est déduite de $A$ par inversion, elle fournit aussi une gamme
 bien répartie, donc à translation près c'est $A$ comme on l'a vu.
 On peut donc énoncer le théorème précédent sous une forme plus simple encore:
 
 \og $A$ s'écrit $\{0, \lfloor \dfrac{c}d\rfloor,… \lfloor \dfrac{(d-1)c}d\rfloor\}$
 à translation près\fg.
\end{Rem}
\begin{Rem}
 Plus généralement, on peut raisonner avec les entiers immédiatement
 supérieurs, ou l'entier le plus proche, au lieu des parties entières. 
 En fait, ainsi que le fait remarquer \cite{CD}, une forme plus générale
 quoique équivalente de la séquence donnant $A$ est
$$k\mapsto J_{c,d}^\alpha (k) =\lfloor \dfrac{k c + \alpha}d \rfloor$$
avec $\alpha$ réel arbitraire. Ce n'est pas difficile de le prouver en
reprenant la discussion précédente adaptée à cette nouvelle famille.
Le choix $\alpha=1/2$, éventuellement combiné avec l'inversion de $A$ en -A (ce qui change de fait
$k$ en $c-k$) permet effectivement de passer de l'une à l'autre des trois façons
usuelles d'approximer un réel par un entier.
\end{Rem}
\subsubsection{Le cas $d\wedge c>1$}
\begin{Thm} \label{nonPremier}
    Soit $A = \{ a_{1}, … a_{d}\}$ une gamme bien répartie.
    Alors $A$ est engendrée par une formule comme ci-dessus:
    $$A = \{ \lfloor \dfrac{k d + \alpha}c \rfloor \mid k=0…d-1\}$$
    Quand $m=c\wedge d $ est strictement supérieur à 1, 
     $A$ est $c/m-$ périodique (un mode à transposition limitée, pour reprendre
    la notion du compositeur Olivier \textsc{Messiaen}), obtenue comme orbite 
    sous la translation de $c/m$ d'un domaine fondamental qui est lui-même une
    gamme bien répartie monogène de cardinal $d/m$ dans $\mathbf Z_{c/m}$.
\end{Thm}
\begin{proof}
   Remarque: il est tentant de faire cette démonstration à coups de structures quotients et
   de suites exactes, mais nous avons préféré des arguments plus concrets, étant donnée la
   nature appliquée des objets considérés.
   
   Seul le cas $m=c\wedge d >1$ reste à considérer, l'autre étant élucidé.
   Il n'est plus possible que tous les éléments de $d\, A$ soient consécutifs dans $\zc$, car
   
   \begin{Lemme}
     L'application $\varphi:x\mapsto d\, x$ de $\zc$ dans lui-même a pour image le sous-groupe
     $\pmb{m\, \zc}\approx \mathbf{Z}_{c/m}$ cyclique d'ordre $c/m$.
     Chaque élément de l'image a $m$ antécédents, translatés les uns des autres
     de $c/m$: $\ker \varphi = c/m\, \zc$.
   \end{Lemme}
   Pour maximiser $IC(A)(d)=|\mc{F}_{A}(d)|^2$, le mieux que l'on puisse obtenir
   est que les éléments de $d\, A$ soient consécutifs \textbf{dans $\pmb{m\, \zc}$}, 
   i.e.
   $$
     d\, A = m\, \{0,1,2 … \dfrac{d}m-1\}
          = m\, \{0,1,2 … d'-1\} \quad\text{en posant }d' = d/m.
   $$
   Attention: $A'=d\, A$ a alors $d'=d/m$ elements, et pas $d$ comme $A$. C'est qu'en fait
   on peut voir $d\, A$, image de $A$ par $\varphi$, non comme un
   ensemble mais comme un \og multiset\fg: chaque élément
   de $d\, A$ est répété $m$ fois, comme image de $m$ éléments distincts de $A$.   
   
   
   Donnons un exemple: la configuration
   $4\times\{0,3,5,8\}=\{0,2,0,2\}\subset\mathbf{Z}_{10}$ 
   est optimalement serrée, plus que $\{0,2,4, 2\}$ ou toute autre configuration plus
   \og lâche\fg, dans $4\,\zd = 2\,\zd$.
   
   On prouve que la configuration $ d\, A = m\, \{0,1,2 … \dfrac{d}m-1\}$ 
   (au sens où chaque élément est répété $m$ fois) par le lemme des points bien serrés: 
   si le multiset $d\, A$ \textbf{n'est pas} constitué de $d'=d/m$ éléments consécutifs répétés
   chacun $m$ fois, alors on peut augmenter la somme des exponentielles
   en rapprochant un des points de la somme de tous. On itère jusqu'à obtenir une
   configuration maximale (cf. figure).
      
   Ceci étant acquis, nous avons vu comment obtenir
   $d' A' = \{0,… d'-1\}$ modulo $c' = c/m$: il suffit de poser
   $$
     A' = -\{0, \lfloor \dfrac{c'}{d'}\rfloor,… \lfloor \dfrac{(d'-1)c'}{d'}\rfloor\}
     =-\{0, \lfloor \dfrac{c}{d}\rfloor,… \lfloor \dfrac{(d'-1)c}{d}\rfloor\}
  $$
  N.B.: $c', d'$ sont premiers entre eux, et $A'$ relève donc du cas précédent: 
  c'est une gamme bien répartie dans $\mathbf Z_{c'}$ (un sous-cycle, \og subcycle\fg, 
   dans la terminologie de \textsc{Cohn}, cf. \cite{Quinn}).
  
  On tient alors $A$ tout entier en rajoutant à $A'$, plongé dans $\zc$, le noyau de $\varphi$,
  i.e. en translatant $A'$ des multiples de $c'=c/m$, car comme on l'a dit 
  $$
    a\in A \iff \exists k\in [0, d'-1],\ d a = k \pmod c\iff 
   d' a = m^{-1} k \pmod c' \iff a\in A' \pmod c'=c/m
  $$
  En résumé,
  $$
      A    = A' \oplus \{0, \dfrac{c}m, \dfrac{2c}m, …\}
   $$
  Donc $A$ est un Mode à Transpositions Limitées\footnote{Pour
   prendre une définition mathématique, et non la définition (musicale !) de \textsc{Messiaen}, 
   nous dirons qu'un M.T.L.
  est une partie  de $\zc$ dont l'orbite sous l'action des translations a moins
  de $c$ éléments distincts, ou encore qu'elle admet un fixateur non trivial}. 
  De plus, nous avons, encore (modulo $c$),
   $$
     A= -\{0, \lfloor \dfrac{c}{d}\rfloor,… \lfloor \dfrac{(d-1)c}{d}\rfloor\}
   $$
  En effet, pour $k≥d'$ on écrit $k = q d' + r, 0≤r<d'$ d'où
  $$
    \lfloor \dfrac{k c}{d}\rfloor = \lfloor \dfrac{(q d' + r)c'}{d'}\rfloor
    = q c' + \lfloor \dfrac{r c}{d}\rfloor \in (A' + q c')
  $$  
  Nous avons établi (à translation près) que les formules du type \cite{CD}
  engendrent toutes les gammes bien réparties.
\end{proof}
\begin{Rem}\label{nombreBGR}
  Le nombre de gammes bien réparties distinctes est $c'=c/m$ (ce n'est $c$ que si
  $d\wedge c = 1$), et toute \textbf{inversion} d'une gamme bien répartie
  est toujours une \textbf{translatée} d'elle même: le groupe des inversions-translations
  (i.e. les applications de la forme $x\mapsto \lambda \pm x$ dans $\zc$)
  n'agit pas fidèlement sur l'ensemble des gammes bien réparties.
  
  Par exemple, le prototype de $GBR_{10,4}$ est $(0,3,5,8)$.
  La gamme rétrogradée est $10 - (0,3,5,8) = (2,5,7,0) = (0,3,5,8) + 2$.
\end{Rem} 
Notons dans le même esprit que la proposition \ref{inclusComplement} reste vraie
pour ces GBR générales (non monogènes): une réduction au quotient modulo $c'$
ramène au cas $c',d'$, et (pour $d'<c'/2$) d'un $GBR(c',d')\subset GBR(c', c'-d')$ on déduit
$GBR(c,d)=GBR(c',d')+c'\Z \subset GBR(c', c'-d') + c'\Z = GBR(c,c-d)$.

\subsection{Une propriété des intervalles consécutifs}\

À partir des formules ci-dessus, on obtient bien, mais comme conséquence, 
 la propriété de \textsc{Myhill}) qui permet de retrouver la définition américaine 
 (\cite{CM}, \cite{CD}) des gammes monogènes \textbf{bien formées}
 (\og Well Formed\fg):

\begin{Thm}
   Les intervalles entre deux notes consécutives d'une gamme bien répartie
   ne peuvent pas prendre plus de deux valeurs distinctes. Il y a une seule valeur
   si et seulement si la gamme divise également $\zc$.
\end{Thm}
Ce dernier cas est considéré comme dégénéré (\og not Well Formed\fg).
Remarquons que ce résultat est valide aussi pour les intervalles entre notes
\textbf{d'indices} à une distance constante, le cas des notes consécutives
étant celui d'une différence d'indices de 1. C'est sous cette forme plus forte qu'est
généra\-le\-ment énoncée la propriété de \textsc{Myhill}. Clarifions par un exemple :
dans le cas de la gamme majeure $GBR_{12,7}$, dont une instance est
$(0,2,4,5,7,9,11) = $(do, ré, mi, fa, sol, la, si),
les intervalles (chromatiques)
entre deux notes consécutives sont de deux demi-tons, sauf mi-fa et si-do
qui sont d'un seul demi-ton: il existe des secondes majeures et des secondes mineures, 
contrairement à d'autres gammes non bien réparties, comme la gamme mineure 
par exemple $(0,2,3,5,7,8,11)$ qui contient aussi des secondes augmentées.
Revenons à la gamme majeure : les intervalles de deux en deux sont les 
\textbf{tierces} do-mi, ré-fa, mi-sol etc…
qui contiennent trois (tierces mineures) ou quatre (tierces majeures) demi-tons.
Il y a de même deux sortes de quartes (do-fa et fa-si), de quintes, etc…
Revenons à la preuve du théorème:

\begin{proof}
  Les deux valeurs possibles sont tout bonnement les entiers immédiatement voisins
  de $c/f$ où $f$ est le générateur dans les formules précédentes.
 En effet, on a par l'encadrement classique de la partie entière
 $$
    \dfrac{c}d - 1 = \dfrac{(k+1)c}{d} - 1 - \dfrac{k\,c}{d} <
    \lfloor \dfrac{(k+1)c}{d}\rfloor - \lfloor \dfrac{k\,c}{d}\rfloor
    < \dfrac{(k+1)c}{d} - (\dfrac{k\,c}{d}-1)
    = \dfrac{c}d + 1
 $$
 et cet intervalle ouvert est de largeur 2, donc contient deux entiers, nommément
 $\lfloor \dfrac{c}d \rfloor$ et  $\lceil \dfrac{c}d \rceil$.
 
 Les deux sont possibles quand $c/d$ est non entier: le plus grand quand
 $\dfrac{k\,c}{d}$ est non entier, le plus petit par exemple pour $k=0$.
 
 Si $c/d = m$ est un entier en revanche, on a toujours 
 $\lfloor \dfrac{k\,c}{d}\rfloor = \dfrac{k\,c}{d}=k\,m$
 donc l'intervalle entre deux éléments consécutifs vaut toujours $m$.
 
 Nous avons établi que les intervalles de \og seconde\fg\ (entre deux
 notes consécutives d'une même gamme) peuvent prendre seulement
 deux valeurs distinctes, qui sont les deux arrondis d'un même rationnel.
 Enfin, on trouve de même les intervalles de \og tierce, de quarte…\fg, en considérant
 les arrondis de $\dfrac{2c}d, \dfrac{3c}d$…
\end{proof}
Par exemple, les intervalles consécutifs entre les éléments de la gamme majeure
($(c,d)=(12,7)$) sont 2 (cinq fois) et 1 (2 fois).

\subsection{Cardinal et variété}\ 

Une autre jolie propriété qui a originellement amené \cite{CM} à leur définition (dans le cas
où $c\wedge d = 1$ seulement) est équivalente à la propriété de \textsc{Myhill}
évoquée ci-dessus (qu'elle généralise), et facile à établir comme conséquence des formules
génératrices retrouvées dans le présent exposé. Donnons là d'abord sur un exemple : 
si on cherche les $GBR_{7,3}$ \textbf{dans l'univers diatonique}, c'est à dire
en numérotant les notes de la gamme de do majeur, on obtient 7 accords, 
do-mi-sol, ré-fa-la, etc…
Ces accords se trouvent répartis en trois types (majeur, mineur, diminué).
Trois notes, trois types: la cardinalité est égale à la variété, i.e. au nombre d'orbites modulo 
transposition. Ce concept est musicalement très important car lié à la notion d'ambiguïté.
Mathématiquement il est sans doute moins naturel que la notion de \og bonne répartition\fg.

\begin{Thm}\cite{CM2}
   Considérons pour $c\wedge d =1$ une $GBR_{c,d}$, que l'on écrit sous la forme
   \textbf{ordonnée} $(a_{1},… a_{d})$. Alors pour tout $k<d$, il y a exactement
   $k$ orbites par translation des $d$ séquences $S_{k,i} =(a_{i}, a_{i+1}, … a_{i+k-1})$ 
   où les indices sont pris modulo $d$.
\end{Thm}


\begin{proof}
  Plutôt que de donner la démonstration dans toute sa généralité (cf. \cite{CD}, thm 1.10)
  je préfère mettre en lumière la compréhension fine du phénomène qu'apporte la
  transformation $A'\ni x\mapsto f x \in A$.
  Il arrive en effet fréquemment, mais pas systématiquement, que l'ordre des éléments
  de $A$ soit le même que l'ordre des éléments de $A'$. C'est le cas par exemple quand
  $d<\!\!< c$, comme pour $c=15, d=4$. Alors comme en fait foi la figure, 
  il est clair qu'il n'y a que $d$ types de séquences d'éléments contigus de $A'$, et donc
  (puisqu'on a supposé que l'homothétie $x\mapsto f\times x\mod c$ conserve cette contiguïté)
  il en est de même pour $A$, ce qui est la propriété à démontrer.
  
  Le cas où les éléments de $A'$ correspondant aux séquences d'éléments contigus de
  $A$ ne le sont pas, comme par exemple la gamme pentatonique
  $(0,2,4,7,9)\subset\zd, f=5$ est moins évident géométriquement, mais guère difficile à
  prouver algébriquement avec les fonctions $J_{c,d}^\alpha$.
\end{proof}
Par ailleurs on trouvera dans \cite{CD} le cas $c\wedge d>1$ (le nombre de cas
distincts est alors $d' = d/(c\wedge d)$).

\section{Classification des gammes bien réparties}
L'essentiel de ce paragraphe repose sur la très belle synthèse de \cite{Quinn},
qui reprend et complète \cite{CD}.

Pour chaque couple  $(c,d)$ nous avons caractérisé une unique gamme
(à translation près) $GBR_{(c,d)}$. Nous allons étudier de plus près la structure fine de 
cet objet remarquable.

\subsection{Les trois types de GBR}\

La discussion porte sur les facteurs communs à $c$ et $d$. On pose encore
$m = d\wedge c$.

\begin{Def}
    On a une GBR de type $I$ quand $m=1$. Cette gamme est alors monogène
    et bien formée (WF) au sens expliqué ci-dessus.
\end{Def}
Rappelons que dans ce cas, on obtient cette gamme (à translation près)
en considérant les multiples de l'inverse multiplicatif (modulo $c$) $f$ de $d$,
ou les multiples de $-f$.

\begin{Def}
    On a une GBR de type $II_{a}$ quand $m=d$ i.e. $d\mid c$. 
    Cette gamme est alors monogène et dégénérée, elle divise de façon parfaitement
    régulière $\zc$.
\end{Def}
Ainsi la septième diminuée $D7 = (0,3,6,9)$ divise $\zd$ en 4. De même 
pour la gamme par tons (M1 de Messiaen) $(0,2,4,6,8,10)$.

\begin{Def}
    On a une GBR de type $II_{b}$ quand $1<m = c-d<d$. 
    On a alors affaire au complémentaire d'une gamme de type $II_{a}$.
\end{Def}
C'est immédiat: le complémentaire est de cardinal $m$, qui divise $c$,
et c'est une GBR comme complémentaire de GBR.

\begin{Def}
    On a une GBR de type $III$ dans le cas restant: $1<m<d, m≠ c-d$.    
\end{Def}
Notons que la classe des GBR de type $III$ est stable par complémentation, comme
les deux autres classes.

Les deux dernières classes réunissent toutes les GBR qui possèdent une \textbf{période interne},
c'est à dire qui sont des MTL (cf. thm \ref{nonPremier}). Elles sont donc toutes obtenues comme
réunion de translatées d'une GBR plus petite.
On peut arguer, avec \textsc{Clampitt} \textsl{et alii}, que les GBR de type $I$ sont fondamentaux, au sens
où ils permettent d'écrire tous les autres: on obtient comme on l'a vu un GBR de type $III$
en découpant $c$ en $m$ parties égales (comme le type $II_{a}$) et en glissant dans chaque 
partie le même GBR de type $I$ correspondant à $d'$ notes parmi $c'=c/m$.


\subsection{Existence de gammes de type III}\ 

Le résultat prouvé dans ce dernier paragraphe est inédit.

On a observé (\cite{Quinn}) qu'il n'existe pas de gammes de type $III$ pour $c=12$,
qui est tout de même le cas de référence pour les musiciens (occidentaux).
Pourtant il en existe pour $c≠12$, par exemple pour $c=18$ on a
$GBR_{(18,8)} = (0,2,4,6,9,11,13,15)$ comme on le voit sur la figure précédente.

Bien sûr, il est impossible d'avoir une GBR de type $III$ quand $c$ est premier, car alors
seul le type $I$ est possible (à part les cas très limites de $GBR_{(c,0)}$ ou $GBR_{(c,c)}$).

Jusqu'en juin 2005 on n'avait pas trouvé de $c$ composite, $c>12$, pour lequel il n'existe pas
de GBR de type $III$. Ceci est général:

\begin{Thm}[Amiot, juin 2005]
    Pour $c>12$ composite, il existe toujours $d$ tel que $GBR_{(c,d)}$ soit de type III.
\end{Thm}
Ma première preuve reposait sur la conjecture de Bertrand, outil bien lourd (avec sa fameuse démonstration 
en 17 lemmes dûe à Pafnouti Tchebitchef !) pour le simple lemme suivant:

\begin{Lemme}
    Pour tout $c>12$ non premier, il existe un diviseur $k$ de $c$ et un nombre premier $p<k-1$
    tel que $p$ ne divise pas $k$.
\end{Lemme}
J'ai depuis trouvé une preuve plus élémentaire de ce lemme:

\begin{proof}
   Notons que ce lemme est faux pour $c=12$, car au maximum $k=6$ et les nombres 
   premiers \textbf{strictement} inférieurs à 5 divisent tous 6.
   
   Prenons $c$ composite supérieur ou égal à 25, les valeurs inférieures sont vérifiées à la main. 
   Soit $k$ le plus grand diviseur strict de $c$. Qu'il soit pair ou impair, on l'écrit
   $k = 2n+1$ ou $k = 2 n+2$. Comme $k≥\sqrt c$, on a $k≥5$ et $n≥2$.
   
\begin{itemize}
\item Premier cas: c est une puissance de 2.
On pose alors $k=c/2, p=3$. Marche dès que $c≥8$.
\item
Second cas: $c$ possède un facteur \textbf{impair} $k >3$
(pas forcément premier).
On prend cette valeur pour $k$, et on pose $p=2$. 
C'est la construction la plus sexy de toutes: par exemple
$ (0, 2, 4, 6, 8, 9, 11, 13, 15, 17)$ quand $c=18$.
Cela marche pour $c≥10$.
\item
Dernier cas: quand $c = $ une puissance de 2, fois 3.
En vérité c'est  quand $c=2\times2\times3$ que le lemme s'avère faux.
Nous tenons donc bien le cas problématique, mais il n'est pas difficile:
dès que $c≥24$ on pose $k=c/2$ et $p=5$, qui conviennent.
\end{itemize}
\medskip
   
   Maintenant le théorème découle de ce que 
   $$j \mapsto \lfloor \dfrac{k j}p\rfloor
   =  \lfloor \dfrac{j c}{c p/k}\rfloor
   $$
   donne une GBR, qui est forcément de type $III$ car elle n'est ni de type $I$ ni de type $II$:
   en effet elle est réunion des translatés de $GBR_{(k,p)}$ qui est de type $I$ 
   (ce n'est \textbf{pas} une division régulière de $\mathbf{Z}_{k}$).
\end{proof}

\subsection{Récursivité des GBR}\ 

La gamme pentatonique de $\zd$ est une GBR dans la gamme à 7 notes de $\zd$.
Pour $(c,d)=(7,5)$ on trouve en effet la GBR $(0,1,2,4,5)$.
Si l'on injecte ces valeurs comme indices dans 
la gamme majeure $GBR_{(12,7)} = (0,2,4,5,7,9,11) $ on trouve
$(0,2,4,7,9)$ qui n'est autre que la gamme pentatonique $GBR_{(12,5)}$ !
(complémentaire, ici, de la précédente).

De même, extraire 3 notes selon le schéma $GBR_{(7,3)}$ de la gamme majeure
$GBR_{(12, 7)}$ on obtient (selon décalage) les différents accords parfaits, 
majeur ou mineur.


Ce phénomène ne se produit pas pour toutes les valeurs, mais il est assez remarquable 
pour mériter d'être mentionné. Il est à rapprocher de la construction des gammes 
par réduites successives d'une même fraction continue, cf. \cite{Clampitt}\footnote{ 
Pour mettre le lecteur en appétit, observons que 3/5, 7/12,… sont des réduites
du DFC de $\log_{2}(3/2)$ qui fonde classiquement la gamme pythagoricienne…}. 
Comme l'obervaient déjà Clough et Douthett dans l'article fondateur \cite{CD},
cela n'est pas si étonnant puisque l'on considère des sous-polygones 
\textsl{aussi réguliers que possible} de polygones \textsl{aussi réguliers que possible}.
Cela leur permettait de définir des \textbf{GBR de seconde espèce}, comme éléments d'indices
$i\in GBR(d,e)$ d'un $GBR_{c,d}$. De l'univers entier à ces objets, on passe par
le chromatique (12 notes), le diatonique (7 notes), et la fonction tonale (3 notes)!
Cela a été joliment formalisé par un modèle que Jack Douthett a présenté l'été dernier 
à un colloque dédié
à la mémoire de John Clough, celui des \og beacons\fg : on dispose une lumière centrale
dans une pièce circulaire, munie d'ouvertures régulièrement réparties, imbriquée 
dans une structure similaire, avec la contrainte peu physique qu'un rayon lumineux 
\textbf{glisse} jusqu'à l'ouverture la plus proche dans le sens des aiguilles d'une montre (!).

Ceci permet de générer les fonctions 
$J:k\mapsto \lfloor \dfrac{k d + \alpha}c \rfloor$ obtenues plus haut, 
et surtout de les imbriquer les unes dans les autres. Les cas où l'on obtient des ensembles
moins bien équilibrés ne manquent pas d'intérêt, puisque l'on obtient ainsi, par exemple,
des séquences d'accords parfaits que l'on retrouve dans la IX\up{e} de Beethoven !

\section{Conclusion}
\subsection{Historique}\ 

Nous présentons ici un historique très sommaire.
À l'origine, \textsc{Clough} et \textsc{Myerson} (\cite{CM}) 
ont inventorié essentiellement le cas $c\wedge d=1$.
Leur démarche est partie de propriétés musicales, comme (Cardinal = Variété) et
la propriété de \textsc{Myhill} (au plus deux valeurs chromatiques différentes pour les intervalles
diatoniques) et aboutissait aux fonctions génératives du type $J$.

C'est avec Jack \textsc{Douthett} que John \textsc{Clough} 
a étendu la propriété \og Maximally Even\fg\ (i.e. \og bien réparti\fg)
à sa forme la plus générale (\cite{CD}). Jack \textsc{Douthett} a par ailleurs démontré que
cette propriété de bonne répartition admettait une infinité de définitions équivalentes,
en termes de fonctions potentielles strictement convexes sur un réseau d'électrons situés sur un cercle
\cite{Ising}. Ainsi de la définition avec les sommes des distances euclidiennes
que nous avions considérée, puis abandonnée\footnote{ Incidemment, ceci
 démontre le théorème de Toth cité en introduction.}; mais aussi par exemple du potentiel
électrostatique d'une répartition d'électrons sur les sites de $\zc$,
avec une interaction coulombienne. En revanche la distance angulaire n'est \textbf{pas}
strictement convexe et l'on comprend (un peu) mieux pourquoi elle ne donne pas les
mêmes ensembles extrémaux.

Ce résultat, qui établit un lien pour le moins surprenant avec la physique des particules,
montre que nous entrons dans une ère nouvelle, où les recherches de théoriciens de la
musique sont susceptibles d'ouvrir de nouveaux horizons dans d'autres sciences.\footnote{
Toutes proportions gardées, il en est de même pour mes propres recherches liant les
canons rythmiques à certains cas de la conjecture de \textsc{Fuglede} \cite{GdM}.}

Ensuite une autre propriété des gammes a été étudiée pour elle-même par Norman
\textsc{Carey} et David \textsc{Clampitt} (\cite{CC}):
la \og Well Formedness\fg \ (gammes bien formées). Dans le contexte de cet article,
il s'agit des GBR dans le cas $c\wedge d = 1$.
Cette notion permet de s'affranchir (plus ou moins) du présupposé d'un total chromatique,
d'un univers ambiant tempéré, mais il s'agit par définition 
du cas particuliers des gammes monogènes.
Une différence ontologique avec les cas précédents est que
l'on y retrouve les gammes \textbf{naturelles}, comme la gamme pythagoricienne engendrée par
le rapport 3/2 et ses puissances\footnote{ Ce qui correspond à placer sur le cercle des
points d'argument multiple de $\log_{2}3/2$}, et le concept de gamme \textbf{bien formée}
permet de considérer comme gammes isomorphes à celle-là toute gamme
monogène de générateur suffisamment voisin\footnote{entre $2^{1/2} $ et $2^{3/5}$
plus précisément.}
de 3/2, comme en fait fois le dessin suivant où le choix du générateur fait que les notes de
la gamme restent dans l'ordre princeps (un polygone \textbf{régulier}, mais étoilé):


On gagne ainsi des propriétés \textbf{topologiques} des gammes (homéomorphisme,
voisinage, mais aussi \og Winding number\fg, i.e. index au sens de Poincaré,
cf. \cite{Cafagna}).
Cela peut se poursuivre par des études de séquences infinies, comme l'est celle des
$n \log_{2} 3/2\mod 1$ qui surgit de la gamme pythagoricienne, mais aussi bien
avec tout irrationnel. Les propriétés de semi-récursivité évoquées plus haut
peuvent s'interpréter par des cheminements dans l'arbre de
\textsc{Stern-Brocot}. Pour cela, avec de belles propriétés d'autosimilarité, voir 
\cite{Noll}.

Enfin et indépendamment de ces développements foisonnants, c'est en recherchant dans les diverses
théories musicales les accords ou gammes les plus \og typiques\fg\ que Ian \textsc{Quinn}
est retombé sur les GBR. Cherchant par ailleurs à définir rigoureusement cette
 typicité en terme de contenu intervallique, il a utilisé les \og balances de Lewin\fg,
 introduites par celui-ci dans son dernier article (\cite{Lewin03}) en un émouvant
 retour aux sources de sa toute première publication. Et c'est Ian \textsc{Quinn}
 qui aura remarqué que les GBR maximisent une valeur d'une transformée de Fourier discrète,
 ce résultat -- le plus récent historiquement -- étant la propriété que le
 présent article a choisi de prendre comme définition des GBR.
 
\subsection{Perspectives}\ 

Il reste à étudier ce que ces outils \textbf{harmoniques}, au sens d'analyse de Fourier, permettent
de préciser dans un contexte plus large, où les intervalles ne
sont plus des diviseurs d'une même totalité (tempérament égal).
Techniquement cela n'est pas trop difficile -- quitte à considérer des distributions et
non plus des fonctions -- mais il est alors nécessaire de réintroduire un peu
 de structure. Cela est possible en prenant
de nouveau en compte un ordre des notes dans la gamme, ce qui fait sens
notamment pour les gammes monogènes (generated scales). On peut
alors définir la transformée de Fourier de la fonction qui à $k\in\Z/d\Z$ associe
la $k\up{e}$ note de la gamme, ce qui n'est \textbf{pas} la définition que nous avons utilisée.

Les recherches en cours ont déjà révélé certaines propriétés géométriques remarquables,
mais pas caractéristiques, des coefficients de Fourier de certaines gammes. 

\subsection{Remerciements}\ 

Merci à D. Clampitt, Thomas Noll, Richard Cohn, Ian Quinn… et à feu David Lewin,
dont les idées continuent de répandre leur lumière.



\end{document}